\documentclass[11pt]{amsart}
\usepackage[margin=1in]{geometry}
 
 \usepackage[margin=1in]{geometry}
\usepackage{hyperref}
\usepackage{amsmath, amscd, amsthm, amssymb}
\usepackage{tikz}
\usepackage{multicol}
\usepackage{framed}
\usepackage{enumitem}
\usetikzlibrary{through,calc}
\usepackage{latexsym}
              
\newtheorem{thm}{Theorem}[section]
\newtheorem{cor}[thm]{Corollary}
\newtheorem{lem}[thm]{Lemma}

\newtheorem{prop}[thm]{Proposition}

\newtheorem*{conj*}{Conjecture}

\newtheorem{defn}[thm]{Definition}
\newtheorem{rmk}[thm]{Remark}

\newcommand{\blem}{\begin{lem}}
\newcommand{\elem}{\end{lem}}
\newcommand{\brop}{\begin{prop}}
\newcommand{\erop}{\end{prop}}
\newcommand{\bcor}{\begin{cor}}
\newcommand{\ecor}{\end{cor}}
\newcommand{\brem}{\begin{rmk}}
\newcommand{\erem}{\end{rmk}}
\newcommand{\benum}{\begin{enumerate}}
\newcommand{\enum}{\end{enumerate}}
\newcommand{\btheo}{\begin{thm}}
\newcommand{\etheo}{\end{thm}}
\newcommand{\bdem}{\begin{proof}}
\newcommand{\edem}{\end{proof}}
\newcommand{\bdefi}{\begin{defn}}
\newcommand{\edefi}{\end{defn}}

\newcommand{\bilif}{\langle \cdot,\cdot \rangle}

\newcommand{\dsl}{\displaystyle \left(}
\newcommand{\rb}{\right)}

\newcommand{\smin}{\smallsetminus}
\newcommand{\sbs}{\subset}

\newcommand{\sps}{\supset}

\newcommand{\R}{\mathbb{R}}

\newcommand{\calC}{\mathcal{C}}
\newcommand{\calF}{\mathcal{F}}
\newcommand{\calP}{\mathcal{P}}

\newcommand{\nin}{\notin}

\newcommand{\rint}{\mathrm{relint}\,}
\newcommand{\relint}{\mathrm{relint}\,}

\title{Periods}
\begin{document}
\renewcommand{\theequation}{\arabic{equation}}
\numberwithin{equation}{section}

\title{A new valuation on polyhedral cones}
\author{Hailun Zheng}
\address{Department of Mathematics\\ University of Michigan\\ Ann Arbor, MI USA}
\email{hailunz@umich.edu}

\author{Micha\l{} Zydor}
\address{Department of Mathematics\\ University of Michigan\\ Ann Arbor, MI USA}\email{zydor@umich.edu}

\begin{abstract} We define a new family of valuations on 
polyhedral cones valued in the space of bounded polyhedra. \end{abstract}
\maketitle

\newcommand{\cpc}{\calP \calC}

\section*{Introduction}

Let $C$ be a polyhedral cone in a $\R$-vector space $V$ and let $\calF(C)$ be the set of all faces in $C$. A classical result in polyhedra geometry, known as the Brianchon-Gram-Sommerville relation (see, for example, \cite{shephard} and \cite{schneider}), reveals the relation between the indicators of the cone $C$ and those of the angle cones $A(F, C)$ in correspondence to all the faces $F$ in $C$:
\[\sum_{F\in \calF(C)}(-1)^{\dim C-\dim F} [A(F,C)]=[-\relint C].\]
There is also a dual version of the above identity: 
\[\sum_{F\in\calF(C)} (-1)^{\dim F-\dim F_0} [F^*]=[-\rint C^*],\]
where $F_0$ is the minimal face of $C$. 

This motivates us to define a function with two variables that evaluates on both the angle cones and the dual cones of the faces. For a cone $C$ we let 
\[\Gamma(C,x, y)=\sum_{F\in\calF(C)}(-1)^{\dim F}[A(F,C)](x)[F^*](y). \]
Furthermore, we define a function $\Gamma_y$ that sends the indicator $[C](x)$ to the function $\Gamma(C, x,x-y)$. Our main theorem is the following.

\btheo[cf. Theorem \ref{thm:main}.]
For any $y\in V$, $\Gamma_y$ is a valuation from the algebra of polyhedral cones to the algebra of bounded polyhedra.
\etheo

If $C$ is a cone in $\R^{2}$ and $y$ belongs to the interior of $C \cap C^{*}$, the function $\Gamma_{y}(C)$ 
is a characteristic function of $C \cap (y - \rint C^{*})$.

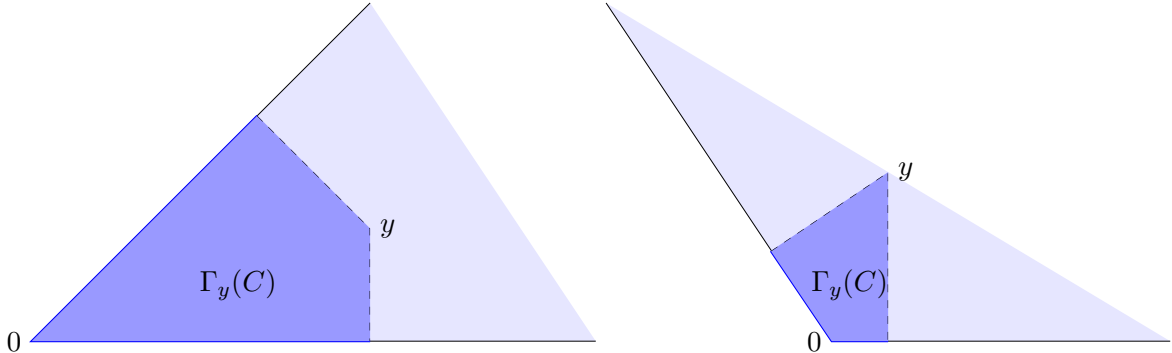
\begin{center}
\begin{tikzpicture}[scale=3]
\draw[very thin, fill=blue!10] (0.5, 1.5) -- (-1, 0) -- (1.5,0);
\draw (0.5,0.5) node[right] {$y$};
\draw (-1,0) node[left] {$0$};
\draw[thick, blue]  (0,1) --  (-1,0) -- (0.5, 0);
\clip  (-1,0) -- (0,1) -- (0.5,0.5) -- (0.5, 0) -- cycle;
\filldraw[color=blue!40] (-1.4,-1.4) rectangle (1.4,1.4);
\draw[dashed, black] (0,1) -- (0.5,0.5)  -- (0.5, 0);
\draw (-0.3,0.25) node[right] {$\Gamma_{y}(C)$};
\end{tikzpicture}
\begin{tikzpicture}[scale=3]
\draw[very thin, fill=blue!10] (-1, 1.5) -- (0, 0) -- (1.5,0);
\draw (0.25,0.75) node[right] {$y$};
\draw (0,0) node[left] {$0$};
\draw[thick, blue]  (-0.27,0.4) --  (0,0) -- (0.25, 0);
\clip  (-0.27,0.4) -- (0,0) -- (0.25, 0)-- (0.25,0.75)-- cycle;
\filldraw[color=blue!40] (-1.4,-1.4) rectangle (1.4,1.4);
\draw[dashed, black] (-0.27,0.4)  -- (0.25,0.75) -- (0.25, 0);
\draw (-0.135,0.25) node[right] {$\Gamma_{y}(C)$};
\end{tikzpicture}
\end{center}

The function $\Gamma_{y}$ (or its close cousin, discussed in Paragraph \ref{ssec:rintC}) 
was first introduced by Arthur \cite{arthur2} (see also Paragraph 1.7 of \cite{labWal}) in his work on the Arthur-Selberg trace formula on locally symmetric spaces. 
Arthur studied it for simplicial cones that arise as chambers associated to parabolic subgroups of 
a reductive algebraic group. 
In this context, $\Gamma_{y}$ allows to decompose the domain of integration (the Siegel domain) into manageable pieces.
Our work can be seen as a generalization of Arthur's construction to arbitrary cones. It will be used in the second author's 
work \cite{zyd} that aims to generalize part of Arthur's work to a more general class of integrals (period integrals).

\section{Preliminaries}\label{sec:prelims}
We begin with introducing basic notations in the algebra of polyhedra. Let $V$ be a fixed finite dimensional Euclidean space over $\R$. Let $\bilif$ be the scalar product on it. For $y \in V \smin \{0\}$, denote by $H_{y}$ the hyperplane in $V$ perpendicular to $y$, and by
\[
H_{y}^{+} = \{v \in V \ | \ \langle v,y \rangle \ge 0\}, 
\quad H_{y}^{-} = -H_{y}^{+}
\]
the closed halfspaces determined by $H_{y}$.
In the following a halfspace always means a set $H^{+}$ of the form $H_{y}^{+}$. A
\emph{cone} (resp.  \emph{polyhedron}) is an intersection of 
a finite number of halfspaces (resp. translates of halfspaces). 
In particular, a cone is always non-empty as $0$ is an element of every cone.

Associate a polyhedron $P$ with the \emph{indicator function} $[P]$. Let $\calP(V)$ be the $\R$-vector space spanned by the indicator functions of all polyhedra $P \subset V$. 
Let $\cpc(V)$ be the vector space spanned by the indicator functions 
of cones. Finally let $\calP_{bd}(V)$ be the vector space spanned by the indicator functions of bounded polyhedra. It turns out that $\calP(V), \cpc(V)$ and $\calP_{bd}(V)$ possess a natural algebra structure, see \cite{Barvinok} for a reference.
Let $W$ be an $\R$-vector space. A linear transformation $T : \calP(V) \to W$ is called a \emph{valuation}. For example, the Euler characteristic $\chi: \calP(V)\to \R$ such that $\chi([P])=1$ for all non-empty polyhedron $P$ is a valuation. The following theorem is a useful tool to decide whether a map on $\cpc(V)$ is a valuation. It was originally stated in \cite[Theorem 2]{Groemer} for polytopes but it continues to hold in the following form for cones.

\btheo\label{thm: valuation IEP}A map $T : \cpc(V) \to W$ is a valuation, if and only if
for all $y \in V \smin \{0\}$ 
and all cones $C$ we have
\[
T(C) = T(C \cap H_{y}^{+}) + T(C \cap H_{y}^{-}) - T(C \cap H_{y}).
\]
\etheo

We will mostly work with cones. 
A \emph{face} of a cone is a cone of the form $C \cap H_{y}$ 
where $C \sbs H_{y}^{+}$. We allow in this definition $y=0$ so that $C$ is always its own face. 
Denote by $\calF(C)$ the set of all faces of $C$. Note that $\calF(C)$ is a singleton $\{C\}$ if and only if $C$ is a subspace of $V$. Let $F_{0} = F_{0}(C)$ be the minimal face of a cone $C$. It is a vector subspace of $V$. 
We call $C$ \emph{pointed} if $F_{0}$ is reduced to the origin. For each face $F$ in $\calF(C)$, let $V_{F}$ be the subspace of $V$ spanned by $F$ and let $d_{F}$ be its dimension. Define $\varepsilon_{F}=(-1)^{d_F}$ and $\varepsilon_{C}^{F}=(-1)^{d_C-d_F}$. We recall the Euler characteristic formula
\begin{equation}\label{eq:euler}
\sum_{F \in \calF(C)} 
\varepsilon_{C}^{F} = 
\begin{cases}
1, & \text{ if }C\text{ is a subspace of }V, \\
0  & \text{ else}. 
\end{cases}
\end{equation}
In the following we will usually use $C$ to denote a cone and use $E, F, G$ to denote its faces.
\subsection{Angle cones and dual cones}

Let $C$ be a cone. Let $\rint C$ be the largest open subset of $V_{C}$ contained in $C$. 
Define the \emph{dual cone} $C^*$ of $C$ as
\[
C^{*} =  \{ v \in V \ | \ \langle v , C \rangle \ge 0\}.
\] 
Also define 
the \emph{angle cone} as
\[
A(F,C) = \{ a(x-z)  \ | \ a > 0, \, x \in C, z\in \rint F \}.
\]
The \emph{Minkowski sum} of two cones $C_1$, $C_2$ is defined to be 
\[
C_1 + C_2 = \{x + y  \ | \  x \in C_1, \, y \in C_2 \}.
\]

We state some well known properties of the angle cone and the dual cone.

\blem\label{lem:wellknown} Let $C$ be a cone, $F_0$ the minimal face of $C$ and $F$ an arbitrary face of $C$.
\begin{enumerate}
\item $A(F_{0},C) =C$, $A(C,C) = V_{C}$ and $(F_0)^*=F_0^\perp$.
\item (\cite[Theorem 4.13]{Barvinok}) $C$ is the Minkowski sum of $F_0$ and some pointed cone $C'$. In particular, $C^*=(F_0+C')^*=F_0^*\cap (C')^*$.
\item There is an inclusion-preserving bijection
\begin{equation*}
	\begin{split}
	 \{\,G\in \calF(C): G\supset F \,\} & \xrightarrow{\sim} \calF(A(F,C))\\
	   G &\mapsto A(F,G).
	\end{split}
\end{equation*} 
Furthermore, for any face $G\in \calF(C)$ containing $F$, we have $A(A(F,G), A(F,C)) = A(G,C)$.
\item 
There is an inclusion-reversing bijection
\begin{equation*}
	\begin{split}
	 \{\,G\in \calF(C): G\supset F \,\} & \xrightarrow{\sim} \calF(A(F,C)^{*})\\
	   G &\mapsto A(G, C)^{*}.
	\end{split}
\end{equation*} 
Moreover, for any face $G \in \calF(C)$ containing $F$, we have $A(A(G,C)^{*}, A(F,C)^{*}) = A(F,G)^{*}$.
\item (\cite[Theorem 5.3]{Barvinok}) There exists a unique valuation $D: \cpc(V)\to \cpc(V)$ such that $D([C])=[C^*]$ for all cones $C$.
\end{enumerate}	
\elem

We have the following classical result, usually referred to as the  Brianchon-Gram-Sommerville relation \cite{shephard} (see also \cite{schneider}, Theorem 3.1)

\btheo\label{thm:bgs}  The following identities hold
\[
\sum_{E \in \calF(C)}\varepsilon_{C}^{E}[A(E,C)] = [- \rint C]
\]
and
\[
\sum_{E \in \calF(C)}\varepsilon_{E}^{F_{0}} [ E^{*}] = [-\rint C^{*}].
\]
We also have for all $F \in \calF(C)$
\[
\sum_{E \sps F}\varepsilon_{C}^{E}[\rint A(E,C)] = [- A(F,C)]
\]
and
\[
\sum_{E \in \calF(C)}\varepsilon_{E}^{F_{0}} [\rint  E^{*}] = [-C^{*}].
\]
\etheo

\section{Main results}

By Theorem \ref{thm:bgs}, both \[\sum_{E \in \calF(C)}\varepsilon_{C}^{E}[A(E,C)] \quad \mathrm{and}\quad \sum_{E \in \calF(C)}\varepsilon_{E}^{F_{0}} [ E^{*}]\] are valuations. This motivates us to extend them to a new function with two variables, called the $\Gamma$ function. We will see in the next subsection that a one parameter family of valuations can be derived from the $\Gamma$ function.

\subsection{The $\Gamma$ function}

\bdefi For any cone $C$ and any points $x, y\in V$, let 
\[
\Gamma(C, x, y) = \sum_{F\in \calF(C)} \varepsilon_{F} [A(F,C)](x)[F^{*}](y).
\]
\edefi

We begin with computing the value of $\Gamma(C, x,y)$ for some special $x,y\in V$. When the context is clear, we abbreviate $\Gamma(C, x,y)$ as $\Gamma(x,y)$.

\blem\label{lem:gamEazy} 
\begin{enumerate}
\item If  $y \in C^{*}$ 
or if $x \in -\rint C$ 
then $\Gamma(x,y) = \varepsilon_{C} [-\rint C](x)[ C^{*}](y)$. 
\item If  $y \in -\rint C^{*}$ 
or if $x \in  C$ then $\Gamma(x,y) = \varepsilon_{F_{0}} [C](x)[-\rint C^{*}](y)$.
\item If $C$ is a vector subspace of $V$, then 
$\Gamma(x,y) = \varepsilon_{C} [C](x) [C^{\perp}](y)$.
\item $\Gamma(x,y) = 0$
unless $x \in V_{C}$  and $y \in F_{0}^{\perp}$.
\end{enumerate}
\elem
\bdem 
The first two identities are immediate consequences of the first two identities in Theorem \ref{thm:bgs} and Lemma \ref{lem:wellknown} (3)-(4). The last two identities follow from Lemma \ref{lem:wellknown} (1)-(2).
\edem

We say that the cone $C$ is \emph{non-degenerate} if it is pointed and $d_C = \dim V_C = \dim V$.

\brop\label{prop:xyvanish} Suppose $C$ is non-degenerate. 
If $\langle x, y \rangle > 0 $ then $\Gamma(x,y) = 0$.
\erop
\bdem 
Suppose $\langle x, y \rangle > 0 $.
We can assume we are not in any of the first three cases considered 
in Lemma \ref{lem:gamEazy}, since in these cases the result is clear. 
Moreover, without changing the value of $\Gamma(x,y)$ 
or the fact that $\langle x, y \rangle > 0 $ we can assume 
that $x$ and $y$ are in general position with respect to $C$, that is 
to say, both $x$ and $y$ do not belong to any space $V_{F}$ 
or $V_F^{\perp}$ for $F \in \calF(C) \smin \{ \{0\},C\}$. 

Let $C_{y} = H_{y} \cap C$ 
and $C_{y}^{+} = H_{y}^{+} \cap C$.
Since $y\notin C^*\cup (-\rint C^*)$, $C_{y}^{+} $ is a non-degenerate cone in $V$ 
and $C_{y}$ is its facet.
Let $\calF(C, y)$ be the subset of $\calF(C)$ 
composed of faces $F$ of dimension at least $2$, such that 
$\dim (F \cap H_{y}) = \dim F - 1 \ge  1$. 
The map 
\[
F \in \calF(C, y) \mapsto F \cap H_{y} 
\]
is a bijection between $\calF(C, y)$ and $\calF(C_y) \smin \{  \{0\} \}$.

We assume that $y$ is in general position with respect to $C$ 
so $[F^{*}](y) \neq 0$ for a face $F$ of $C$ means that $F\in \calF(C_{y}^{+})$ 
and $F \nin \calF(C, y)$. Hence
\[\Gamma(C, x, y) = \sum_{F\in \calF(C)} \varepsilon_{F} [A(F,C)](x)[F^{*}](y)=\sum_{F\in \calF(C_{y}^+)\smin \calF(C,y)}\varepsilon_{F}[A(F, C_{y}^+)](x).\]
The rest of the faces in $C_y^+$ have nontrivial intersection with $H_y$. So they are either $\{0\}$, or those in bijection with faces in $\calF(C,y)$. Applying the first identity of Theorem \ref{thm:bgs} to $C_{y}^{+}$, we obtain that 
\begin{equation*}
\begin{split}
\varepsilon_{C}[-\rint C_{y}^{+}](x) &= \sum_{F\in \calF(C_y^+)}\varepsilon_{F}[A(F, C_y^+)](x)\\
&=[C^{+}_{y}](x)+ \sum_{F \in \calF(C, y)} \varepsilon_{F} \left(
[A(F\cap C_y^+, C_{y}^{+})](x) - [A(F \cap H_{y}, C_{y}^{+})](x)\right) +\Gamma(C,x,y). \\ 
\end{split}
\end{equation*}
The fact that $x \nin -\rint C$ implies that $[-\rint C_{y}^{+}](x) = 0$ and the fact that $x \nin C$ implies that 
$[C_{y}^{+}](x) = 0$. 
We obtain hence
\[
\Gamma(C,x,y) =  \sum_{F \in \calF(C, y)} \varepsilon_{F}\left( [A(F \cap H_{y}, C_{y}^{+})](x) - [A(F\cap C_y^+, C_{y}^{+})](x)\right) .
\]
Since $\langle x, y \rangle > 0$, it follows that $[A(C_{y}, C_{y}^{+})](x) = [A(C_{y}^+, C_{y}^{+})](x)=1$. For all other faces $F \in \calF(C, y)\smin C$, we have 
\[
[A(F \cap H_{y}, C_{y}^{+})](x) = [A(F \cap C_{y}^{+}, C_{y}^{+})](x),
\]
which implies that $\Gamma(C, x,y) = 0$.
\edem

\subsection{The one parameter family of valuations}

\bdefi Let $C$ be a cone in $V$. For $y \in V$, let 
$\Gamma_{y}([C])(x)= \Gamma(C,x,x-y)$ 
where $x \in V$. 
\edefi
First we show that $\Gamma_y([C])$ is supported in a disk of radius $\frac{\|y\|}{2}$ in $V$.  For $v \in V$ and $r \in \R_{\ge 0}$ 
let $B(v,r) = \{w \in V \ | \ \|w-v\| \le r\}$.

\blem\label{lem:xyvanish} 
$\Gamma_{y}([C])(x)  = 0$ unless $\langle x, x-y \rangle \le 0$.
\elem
\bdem 

By Lemma \ref{lem:gamEazy} (2), there exists a pointed cone $C'$ such that $C=F_0+C'$, where $C'$ is in the orthogonal complement of the subspace $F_0$. Furthermore, every face $F$ of $C$ can be written as $F_0+F'$, where $F'$ is a face in $C'$. Let $x'$ and $y'$ be orthogonal projections of $x$ and $y$ 
onto $F_0^\perp$. We have that $[A(F_0+F', F_0+C')](x)=1$ if and only if $[A(F', C')](x')=1$. 
Similarly, $[(F_0+F')^*](y)=[F_0^\perp \cap (F')^*](y)=1$ if and only if $[F'^*](y')=1$. Hence 
\[
\Gamma_y([C])(x)=\sum_{F'\in \calF(C')} \varepsilon_{F_0+F'} [A(F_0+F', F_0+C')](x)[(F_0+F')^*](x-y)=\varepsilon_{F_{0}}\Gamma_{y'}([C'])(x').
\]

Since $\langle x, x-y \rangle = 
\langle x', x'-y' \rangle$, it suffices to prove the claim for pointed cones. This immediately follows from Lemma \ref{lem:gamEazy} (4) and Proposition \ref{prop:xyvanish}.
\edem
 
\bcor\label{cor:gamBound} 
$\Gamma_{y}([C])(x)$
 vanishes unless 
 $x \in B\dsl\dfrac{y}{2}, \dfrac{\|y\|}{2}\rb$. 
 \ecor
 \bdem 
Suppose  $x \nin B\dsl\dfrac{y}{2}, \dfrac{\|y\|}{2}\rb$. 
 That means $\|x-y/2\| > 1/2 \|y\|$. Expanding it we get 
 \[
 \|x\|^{2} - \langle x,y \rangle = \langle x, x-y \rangle > 0
 \]
 which by  Lemma \ref{lem:xyvanish} implies $\Gamma_{y}([C])(x) = 0$.
 \edem

We are ready to prove our main theorem.

\btheo\label{thm:main} 
Let $y \in V$. $\Gamma_{y}$ extends to a valuation on  $\cpc(V)$ with values in $\calP_{bd}(V)$.
\etheo
\bdem
The corollary \ref{cor:gamBound} implies that $\Gamma_{y}([C]) \in \calP_{bd}(V)$ for all cones $C$. 

To show that $\Gamma_{y}$ is a valuation, by Theorem \ref{thm: valuation IEP} it is enough to show that for all cones $C$ and all 
hyperplanes $H$ and the corresponding halfspaces $H^{+}$ 
and $H^{-}$ we have
\[
\Gamma_{y}([C])  = \Gamma_{y}([C \cap H^{+}]) + \Gamma_{y}([C \cap H^{-}]) - \Gamma_{y}([C \cap H]).
\]

Fix a hyperplane $H$. 
For all $F \in \calF(C)$ set 
\[
F^{+} = F \cap H^{+}, \quad F^{-} = F \cap H^{-}, \quad F_{H} = F \cap H.
\]
We partition $\calF(C)$ into the following four disjoint subsets.
\begin{enumerate}
\item $\calF_1(C) = \{F \in \calF(C) \ | \ F \cap \rint H^{+} \neq \varnothing, \ F \cap \rint H^{-} \neq \varnothing \}$. 
Note then that $\calF_1(C)$ can be defined alternatively as
\[
\calF_1(C) = \{F \in \calF(C) \ | \ F^{+} \neq F_{H}, \ F^{-} \neq F_{H}\}.
\]
In particular if $F \in \calF_1(C)$, we have $\dim F = \dim F^{+} = \dim F^{-} = 1 + \dim F_{H}$. 
\item $\calF_2(C) = \{F \in \calF(C) \smin \calF_1(C) \ | \  F \cap \rint H^{+} \neq \varnothing\}$.
Note then that $\calF_2(C)$ can be defined alternatively as
\[
\calF_2(C) = \{F \in \calF(C) \ | \ F^{+} \neq F_{H}, \ F^{+} = F\} = 
\{F \in \calF(C) \ | \ F^{+} \neq F_{H}, \ F^{-} = F_{H}\}.
\]
\item $\calF_3(C) = \{F \in \calF(C) \smin \calF_1(C) \ | \  F \cap \rint H^{-} \neq \varnothing\}$.
Note then that $\calF_3(C)$ can be defined alternatively as
\[
\calF_3(C) = \{F \in \calF(C) \ | \ F^{-} \neq F_{H}, \ F^{-} = F\} = 
\{F \in \calF(C) \ | \ F^{-} \neq F_{H}, \ F^{+} = F_{H}\}.
\]
\item $\calF_4(C) = \{F \in \calF(C) \ | \ F = F_{H}\}$.
\end{enumerate}

We identify $C_{H}$ with a face in the cones $C^+$ and $C^{-}$. 
We have then the following identities
\begin{itemize}[label={--}]
\item 
\[
\calF(C^{+}) = \bigsqcup_{F \in \calF_1(C)}\{F^{+}\} \sqcup 
\bigsqcup_{F \in \calF_2(C)}\{F^{+}\} \bigsqcup \calF(C_{H}).
\]
\item 
\[
\calF(C^{-}) = \bigsqcup_{F \in \calF_1(C)}\{F^{-}\} \sqcup 
\bigsqcup_{F \in \calF_3(C)}\{F^{-}\} \bigsqcup \calF(C_{H}).
\]
\item 
\[
\calF(C_H) = \bigsqcup_{F \in \calF_1(C)}\{F_H\} \sqcup 
\bigsqcup_{F \in \calF_4(C)}\{F\}.
\]
\end{itemize}

Let $F \in \calF_1(C)$. By Lemma \ref{lem:wellknown} (5) and the fact that $[F]=[F^+]+[F^-]-[F_H]$, we have
\[
[F^{*}] = [(F^{+})^{*}] + [(F^{-})^{*}] - [F_H^{*}]
\]
This shows that $[A(F,C)](x)[F^{*}](x-y)$ equals
\[
[A(F,C)](x)[(F^{+})^{*}](x-y) + [A(F,C)](x)[(F^{-})^{*}](x-y) - [A(F,C)](x)[F_H^{*}](x-y).
\]
Moreover, $[A(F,C)] = [A(F^{+},C^{+})] = [A(F^{-},C^{-})]$, and it's not hard to see that
\[
A(F,C) \cap H^{+} = A(F_H,C^{+}), \quad 
A(F,C) \cap H^{-} = A(F_H,C^{-}), \quad 
A(F,C) \cap H = A(F_H,C_H).
\]
Hence
\[
[A(F,C)] = [A(F_H,C^{+})] + [A(F_H,C^{-})] - [A(F_H,C_H)].
\]

On the other hand if $F \in \calF_2(C)$ we clearly have
\[
[A(F,C)](x)[F^{*}](x-y) = [A(F^{+},C^{+})](x)[(F^{+})^{*}](x-y),
\]
and if $F \in \calF_3(C)$ we have
\[
[A(F,C)](x)[F^{*}](x-y) = [A(F^{-},C^{-})](x)[(F^{-})^{*}](x-y).
\]

These properties show that the difference of $$\Gamma_{y}([C])(x)$$ 
with $$\Gamma_{y}([C^{+}])(x) + \Gamma_{y}([C^{-}])(x) - \Gamma_{y}([C_{H}])(x)$$
equals the sum over $F \in \calF_4(C)$ of $\varepsilon_{F}[F^{*}](x-y)$ times
\begin{equation}\label{eq:lastone}
[A(F,C)] - [A(F,C^{+})] - [A(F,C^{-})] + [A(F,C_{H})]
\end{equation}
evaluated at $x$. 
It remains to observe that for $F \in \calF_4(C)$ we have
\[
A(F,C) \cap H^{+} = A(F,C^{+}), \quad 
A(F,C) \cap H^{-} = A(F,C^{-}), \quad 
A(F,C) \cap H = A(F,C_H)
\]
which shows that \eqref{eq:lastone} is zero in this case.
\edem

\subsection{A variation: from $[C]$ to $[\rint C]$}\label{ssec:rintC}

Besides the indicator function on cones, sometimes it is more convenient to work with indicator functions on open cones. As mentioned in the introduction, the motivation 
for the valuation $\Gamma_{y}$ comes from the construction of 
Arthur \cite{arthur2}.
Generalizing the definition of $\Gamma$ introduced by Arthur  
faithfully we should have rather studied the sum
\[
\sum_{F}\varepsilon_{F}[\rint A(F,C)](x)[\rint F^{*}](x-y), \quad\mathrm{for}\;\; x,y \in V.
\]
There is a close connection of this expression to the valuation $\Gamma_{y}$. 
To define a variant of the valuation $\Gamma_y$ for open cones, we need the following lemma, that follows easily 
from Theorem \ref{thm: valuation IEP}.

\blem\label{lem:newFromOld} Let $\phi : \cpc(V) \to W$ be a valuation, where $W$ is a vector space. 
Then, $\phi' :\cpc(V) \to W $ defined as
\[
\phi'([C]) := \varepsilon_{C}\phi ([\rint C])
\]
is a valuation.

\elem

By the above lemma, the map
\begin{equation*}
	\begin{split}
		\Gamma'_{y}: \cpc(V) &\to \calP_{bd}(V)\\
		              [C] &\mapsto \varepsilon_{C} \Gamma_{y}([\rint C])
	\end{split}
\end{equation*}
is a valuation. The following can be considered as a reciprocity formula for $\Gamma_y'$.

\brop  Let $x,y \in V$.  Then
\[
\Gamma'_{y}([-\rint C])(x) = \sum_{F\in \calF(C)}\varepsilon_{F}^{F_{0}(C)}[\rint A(F,C)](x)[\rint F^{*}](x-y).
\]
\erop
\bdem
Since $[C]=\sum_{F\in\calF(C)} [\rint F]$ and $[\rint C]=\sum_{F\in\calF(C)}\varepsilon_{C}^F [F]$, using the fact that $\Gamma_y'$ is a valuation we obtain that
\begin{equation*}
\begin{split}
\Gamma'_{y}([-\rint C])&= 
\sum_{F}\varepsilon_{C}^{F}\Gamma_{y}'([-F]) = 
\varepsilon_{C}\sum_{F}\Gamma_{y}([-\rint F]) \\
&= 
\varepsilon_{C}\Gamma_{y}([- C])=\sum_{F}\varepsilon_{C}^{F}[-A(F,C)] [-F^{*}](\cdot - y)\\
&\stackrel{(1)}{=}
\sum_{F}
\varepsilon_{C}^{F}\dsl \sum_{G \sps F}\varepsilon_{C}^{G}[\rint A(G,C)] \rb
\dsl \sum_{E \sbs F} \varepsilon_{E}^{F_0}[\rint (E^{*})](\cdot - y) \rb \\
&=\sum_{E \sbs G}\varepsilon_{C}^{G}\varepsilon_{E}^{F_0}
[\rint A(G,C)] [\rint E^{*}](\cdot - y) \sum_{E \sbs F \sbs G}\varepsilon_{C}^{F}\\
&\stackrel{(2)}{=} 
\sum_{F}\varepsilon_{F}^{F_{0}}[\rint A(F,C)](x)[\rint F^{*}](x-y).
\end{split}
\end{equation*}
Here equation (1) follows from the third and fourth identities in Theorem \ref{thm:bgs} 
as well as Lemma \ref{lem:wellknown} (3)-(4). The equation (2) follows from the fact that $\sum_{E\subset F\subset G} \varepsilon_{C}^F$ is zero unless $E=G=F$ by the Euler characteristic formula \eqref{eq:euler}.
\edem

We end by discussing the functions $\Gamma_0([C])$ and $\Gamma'_0([C])$, i.e. when $y=0$.

\brop\label{prop:langLemma} We have the following identities for all $x \in V$
\begin{enumerate}
\item
\[
\sum_{F \in \calF(C)} 
\varepsilon_{C}^{F} [A(F,C)](x)  [F^{*}](x) = 
\begin{cases}
1, & \text{ if }C\text{ is a subspace of }V\text{ and }x=0, \\
0 & \text{ else}. 
\end{cases}
\]
\item
\[
\sum_{F \in \calF(C)} 
\varepsilon_{C}^{F} [\rint A(F,C)](x)  [\rint F^{*}](x) = 
\begin{cases}
1, & \text{ if }C\text{ is a subspace of }V\text{ and }x=0, \\
0 & \text{ else}. 
\end{cases}
\]
\item
\[
\sum_{F \in \calF(C)} 
\varepsilon_{C}^{F} [A(F,C)^{*}](x^{F}) [F](x_{F}) = 
\begin{cases}
1 , & \text{ if }C\text{ is a subspace of }V, \\
0   & \text{ else}. 
\end{cases}
\]
\item
\[
\sum_{F \in \calF(C)} 
\varepsilon_{C}^{F} [\rint A(F,C)^{*}](x^{F}) [\rint F](x_{F})= 
\begin{cases}
1 , & \text{ if }C\text{ is a subspace of }V, \\
0   & \text{ else}. 
\end{cases}
\]
\end{enumerate}
\erop
\bdem
The Corollary \ref{cor:gamBound}, with $y=0$, immediately yields the first two identities given the Euler relation \eqref{eq:euler}. 
Applying to the first identity 
the valuation that takes characteristic functions of cones to characteristic functions of their duals 
yields then the third identity. The last one is finally an easy consequence of the third one.

\edem

\brem 
The above Proposition, or an essential part thereof, is proved by Schneider \cite{schneider}, Theorem 1.1 
(later proved for general polyhedra in \cite{hugKab}). We provide thus an alternative proof of Schneider's result.
\erem

\bibliographystyle{alpha}

\end{document}